\documentclass[11pt, reqno]{amsart}
\usepackage{mathrsfs}
\usepackage[active]{srcltx}
\usepackage{mathrsfs}
\usepackage{array}
\usepackage{multirow}
\usepackage{tikz}
\usepackage{pgfplots}
\usepackage{soul}
\usepackage{mathtools}
\usepackage{xcolor}
\usepackage{caption}
\usepackage{subcaption}
\usepackage{float}
\usepackage{caption}
\usepackage{yfonts}
\usepackage{color}
\usepackage[unicode]{hyperref}

\numberwithin{equation}{section}

\usepackage{graphicx}

\usepackage[left=2.5cm, right=2.5cm, top=2.5cm, bottom=2.5cm]{geometry}

\newtheorem{theorem}{Theorem}

\newtheorem{claim}[theorem]{Claim}
\newtheorem{lemma}[theorem]{Lemma}

\theoremstyle{remark}
\newtheorem*{remark}{Remark}

\theoremstyle{definition}

\begin{document}
	\title{On possible generalisations of quasi-contractions}
	
	\author{T\"unde Cseh}
	\address{Department of Mathematics, Babe\c s-Bolyai University, Cluj-Napoca, Romania    
	}
	\email{cseh.tunde.0914@gmail.com}

	\author{S\'andor Kaj\'ant\'o}
	\address{Department of Mathematics, Babe\c s-Bolyai University, Cluj-Napoca, Romania  
	}
	\email{kajanto.sandor@math.ubbcluj.ro}
	
	\author{Andor Luk\'acs}
	\address{Department of Mathematics, Babe\c s-Bolyai University, Cluj-Napoca, Romania    
	}
	\email{lukacs.andor@math.ubbcluj.ro}
	
	\date{2021 March 23}
	
	\begin{abstract}
		This paper investigates whether some fixed point theorems for quasi-contractions on metric spaces introduced by \`Cir\`ic in \cite{ciric1974generalization} and generalised by Kumam et al. in \cite{kumam2015generalization} can be improved further. It turns out that the answer is negative. We provide two examples of complete metric spaces and two operators without fixed points. We prove that for any possible straightforward relaxation of generalised quasi-contractive conditions, one of these operators satisfies the condition.	
	\end{abstract}

	\maketitle
	\section{Introduction and preliminary results} 
	Banach's contraction principle is a fundamental result in the study of fixed points of operators defined on complete metric spaces. This principle can be stated as follows.

	Let $(X,d)$ be a complete metric space and $T\colon X\to X$ be a self-map. If there exists $q\in[0,1)$ such that for all $x,y\in X$
	\[d(Tx,Ty)\le q\cdot d(x,y)\tag{C1}\label{eq:c1},\]
	then $T$ has a unique fixed point $x^*\in X$. Furthermore, for any $x_0\in X$ the sequence $x_{n+1}=Tx_n$ converges to $x^*$ in $X$.

	Due to its wide range of applicability in different fields of mathematics, several generalisations have appeared. This paper focuses on one possible ``branch'' of these improvements, that of quasi-contractions, that we present below.

	\subsection{Quasi-contractions}
	The notion was first introduced by \`Cir\`ic in \cite{ciric1974generalization}, hence it is sometimes referred to as \`Cir\`ic-type contractions. It consists of two separate improvements of Banach's original principle. 

	On the one hand, the requirement of completeness of $(X,d)$ is relaxed to $T$-orbitally completeness. We recall that the \emph{orbit} of $T\colon X\to X$ is defined as
	\[O_T(x)=\{x,Tx,\dots, T^n x,\dots\},\]
	and a metric space is \emph{$T$-orbitally complete} if every Cauchy sequence in $O_T(x)$ is convergent in $X$.

	On the other hand, the contractive condition \eqref{eq:c1} is relaxed as well and it is replaced with the following:
	\[d(Tx,Ty)\le q\cdot\max\{d(x,y),d(Tx,x),d(Tx,y),d(Ty,x),d(Ty,y)\}.
	\tag{C2}\label{eq:c2}\]
	With these improvements, the operator $T$ still has a unique fixed point $x^*\in X$ and for any $x_0\in X$ the sequence $x_{n+1}=Tx_n$ converges to $x^*$.

	\subsection{Generalised quasi-contractions} \`Cir\`ic's idea was developed further by Kumam et al. in \cite{kumam2015generalization}. The authors introduced the notion of generalised quasi-contraction which uses the condition
	\begin{align*}
	d(Tx,Ty)&\le q\cdot\max\{d(x,y),d(Tx,x),d(Tx,y),d(Ty,x),d(Ty,y),\\
	&\qquad d(T^2x,x),d(T^2x,y),d(T^2x,Tx),d(T^2x,Ty)\},\tag{C3}\label{eq:c3}
	\end{align*}
	and proved a fixed point theorem similar to \`Cir\`ic's.
	\newpage

	Focusing on conditions (C1-C3), the following questions arise naturally:
	\begin{itemize}
		\item[-] Why is \eqref{eq:c3} not symmetric in $x$ and $y$, i.e. why are the $d(T^2y,\cdot)$ terms excluded? More generally, can one include terms of the form $d(T^ky,\cdot)$, where $k\ge 2$?
		\item[-] Can one introduce additional terms of the form $d(T^kx,\cdot)$, where $k\ge 3$, in the set on the right-hand side?
		\item[-] What is the ``most general'' version of these types of conditions that guarantees the existence and uniqueness of the fixed point of the operator in concern?  
	\end{itemize}
	In the next section we answer all the questions above: the conclusion is that condition \eqref{eq:c3} cannot be relaxed further. 

	\section{Main Result}

	\begin{theorem} 
		There exists a complete metric space $(X,d)$ and an operator $T\colon X\to X$  such that $T$ has no fixed points, while for some $q\in(0,1)$ and for every $x,y\in X$ we have 
		\begin{align*}
		d(Tx,Ty)&\le q\cdot\max\{d(x,y),d(Tx,x),d(Tx,y),d(Ty,x),d(Ty,y),\\
		&\qquad d(T^2x,x),d(T^2x,y),d(T^2x,Tx),d(T^2x,Ty),D\},\tag{C}\label{eq:cx}
		\end{align*}
		where $D$ is one of the  distances
		\[\begin{cases}
			d(T^ax,T^by), \mbox{ for some } a\ge3,b\ge0,\\
			d(T^ax,T^by), \mbox{ for some } a\ge0,b\ge2,\\
			d(T^ax,T^bx), \mbox{ for some } a\ge3, b\ge0, a\ne b,\\
			d(T^ay,T^by), \mbox{ for some } a\ge2, b\ge0, a\ne b.\\
		\end{cases}\tag{$\mathcal{D}$}\label{eq:p}\]
		\label{thm:main}
	\end{theorem}
	The proof is obtained by constructing two different examples, depending on the type of the distance $D$. First, we construct a space and an operator that deals with distances of type $D=d(T^{k+1}x,T^ky)$, for some $k\ge 2$. Then we give another construction that discusses the remaining cases for $D$. We claim that both examples satisfy condition \eqref{eq:cx} with the respective $D$.

		\begin{claim}
		Let $X=\{2^n:n\in \mathbb{N}\}$, $d(x,y)=|x-y|$ and $T(x)=2x$. Obviously $(X,d)$ is a complete metric space and $T$ does not have fixed points. Furthermore, let $D=d(T^{k+1}x,T^ky)$, with $k\ge 2$ arbitrary. Then condition \eqref{eq:cx} holds for all $x,y\in X$. 
		\label{cl:1}
	\end{claim}
	\begin{proof}
		For every $x,y\in X$, there exists $m,n\in\mathbb{N}$, such that $x=2^m$ and $y=2^n$. We have three cases.
	\begin{itemize}
		\item If $m>n$, then 
		\[
			d(Tx,Ty)=2^{m+1}-2^{n+1}\le2^{m+1}-2^{n}= \frac{1}{2} (2^{m+2}-2^{n+1})=\frac{1}{2}d(T^2x,Ty).
		\]
		\item If $m=n-1$, then
		\[
			d(Tx,Ty)=2^{n+1}-2^n=\frac{2}{3}(2^2-1)2^{n-1}=\frac{2}{3}(2^{n+1}-2^{n-1})= \frac{2}{3} d(Ty,x).
		\]
	 	\item If $m<n-1$, then $0\le2^{n-m-1}-2$. Adding $2\cdot 2^{n-m-1}-1$ to both sides, we obtain
	 	\[(2^{n-m}-1)\le 3(2^{n-m-1}-1).\] 
	 	Now using that $k\ge2$, we can write
	 	\begin{align*}
	 		d(Tx,Ty)&=2^{n+1}-2^{m+1}=2^{m+1}(2^{n-m}-1)\le 2^{m+1}\cdot3(2^{n-m-1}-1)\\
	 		&=\frac{3}{2^k}(2^{n+k}-2^{m+1+k})\le\frac{3}{4}(2^{n+k}-2^{m+1+k})= \frac{3}{4}d(T^ky,T^{k+1}x).
	\end{align*}
	\end{itemize}
	 In conclusion, \eqref{eq:cx} holds with $q = \frac{3}{4}$.
	\end{proof}

	\begin{claim}
		 Let $X=\{z^n\mid n\in\mathbb{N}\}$, where $z=-1+i\sqrt{3}$, $d(x,y)=|x-y|$ and $T(x)=zx$. Obviously $(X,d)$ is a complete metric space and $T$ does not have fixed points. Furthermore, let $D$ be one of the distances from \eqref{eq:p}, which is not included in Claim \ref{cl:1}. Then condition \eqref{eq:cx} holds for all $x,y\in X$. 
		\label{cl:2}	
	\end{claim}
	We present two lemmas that we use in the proof of Claim \ref{cl:2}. In the forthcoming proofs, we use the following facts without mention: $|z|=2$, $|z-1|=\sqrt{7}$, $|z^2-1|=\sqrt{21}$, $|z^3-1|=7$ and $|z+1|=\sqrt{3}$.
	\begin{lemma}
		If $z=-1+i\sqrt{3}$, then $D=|z^{u+2}-z^v|\ge \sqrt{21}$, for all $u,v\ge0$, with $u+2\ne v$. \label{lem:1}
	\end{lemma}
	\begin{proof}
	We have the following cases.
	\begin{itemize}
		\item If $u=0$ and $v=0$, then $D=|z^2-1|=\sqrt{21}$.
		\item If $u=0$ and $v=1$, then $D=|z^2-z|=2|z-1|=2\sqrt{7}>\sqrt{21}$.
		\item If $u=0$ and $v=3$, then $D=|z^2-z^3|=4|z-1|=4\sqrt{7}>\sqrt{21}$.
		\item If $u=1$ and $v=0$, then $D=|z^3-1|=7>\sqrt{21}$.
		\item If $u=1$ and $v=1$, then $D=|z^3-z|=2\sqrt{21}>\sqrt{21}$.
		\item If $u=1$ and $v=2$, then $D=|z^3-z^2|=4\sqrt{7}>\sqrt{21}$.
		\item If $u>1$ or $v>3$, then $D=|z^{u+2}-z^v|\ge||z^{u+2}-|z^v||=|2^{u+2}-2^v|\ge8>21$.
	\end{itemize}
	\end{proof}
		\begin{lemma} If $z=-1+i\sqrt{3}$, then $D=|z^{u+3}-z^v|\ge 7$, for all $u,v\ge0$, with $u+3\ne v$. \label{lem:2}
	\end{lemma}
	\begin{proof}
	We have the following cases.
	\begin{itemize}
		\item If $u=0$ and $v=0$, then $D=|z^3-1|=7$.
		\item If $u=0$ and $v=1$, then $D=|z^3-z|=2|z^2-1|=2\sqrt{21}>7$.
		\item If $u=0$ and $v=2$, then $D=|z^3-z^2|=4|z-1|=4\sqrt{7}>7$.
		\item If $u>0$ or $v>3$, then 
		$D=|z^{u+3}-z^v|\ge||z^{u+3}-|z^v||=|2^{u+3}-2^v|\ge8.$
	\end{itemize}
	\end{proof}
	
	\begin{proof}[Proof of Claim \ref{cl:2}] We have four cases.
	\begin{itemize}
		\item If $m=n+s$, with $s\ge 2$ then $d(Tx,Ty)\le q_1d(Tx,x)$, where $q_1=\frac{5}{2\sqrt{7}}<1$. Indeed, we have
		\begin{align*}
		d(Tx,Ty)&=|z^{n+s+1}-z^{n+1}|=2^{n+1}|z^s-1|\le 2^{n+1}(|z^s|+1)=2^{n+1}(2^s+1)\\
		&\le2^{n+1}(2^s+2^{s-2})=5\cdot2^{n-1+s}=\frac{5}{2\sqrt{7}}\cdot\sqrt{7}\cdot2^{n+s}\\
		&=\frac{5}{2\sqrt{7}}|z-1||z^{n+s}|=\frac{5}{2\sqrt{7}}d(Tx,x).
		\end{align*}
		\item If $n=m+s$, with $s\ge 2$ one can similarly prove that 
		$d(Tx,Ty)\le \frac{5}{2\sqrt{7}} d(Ty,y)$.
		\item If $m=n+1$, then $d(Tx,Ty)\le q_2d(T^2x,Ty)$, where $q_2=\frac{\sqrt{3}}{3}<1$. Indeed, we have
		\begin{align*}
		d(Tx,Ty)&=|z^{n+2}-z^{n+1}|=2^{n+1}|z-1|=2^{n+1}\sqrt7\\
		&=\frac{\sqrt3}{3}\cdot\sqrt{21}\cdot 2^{n+1}=\frac{\sqrt3}{3}|z^2-1||z^{n+1}|=\frac{\sqrt3}{3}d(T^2x,Ty).
		\end{align*}
		\item If $m=n-1$, then there exists $q_3\in(0,1)$, such that $d(Tx,Ty)\le q_1D$ and $D$ is any distance from \eqref{eq:p} that was not considered in Claim \ref{cl:1}. 

		To prove this statement, we observe that $D$ can have the following forms.
		\begin{itemize}
			\item If $D = d(T^ax, T^by)$ for some $a\ge 3, b\ge 0$, $a\ne b+1$,  then $D=|z^{n-1+a}- z^{n+b}|$.
			\item If $D = d(T^ax, T^by)$ for some $a\ge 0$, $b\ge 2$, $a\ne b+1$, then $D=|z^{n-1+a}- z^{n+b}|$.
			\item If $D = d(T^ax, T^bx)$ for some $a\ge 3, b\ge 0$, $a\ne b$, then $D=|z^{n-1+a}- z^{n-1+b}|$.
			\item If $D = d(T^ay, T^by)$ for some $a\ge 2, b\ge 0$, $a\ne b$, then $D=|z^{n+a}- z^{n+b}|$.
		\end{itemize}
		This implies that
		\begin{itemize}
			\item either $D = |z^{n+a}- z^{n+b}|$ for some $a\ge 2, b\ge 0$, $a\ne b$,
			\item or $D = |z^{n+a}- z^{n-1+b}|$ for some $a\ge 2, b\ge 0$, $a+1\ne b$.
		\end{itemize}
		On the one hand, using Lemma \ref{lem:1} we have
		\begin{align*}
		d(Tx,Ty)&=|z^{n}-z^{n+1}|=2^n\sqrt7=\frac{\sqrt{3}}{3}2^n\sqrt{21}\\
		&\le\frac{\sqrt{3}}{3}|z^{n}||z^{a}-z^{b}|=\frac{\sqrt{3}}{3}|z^{n+a}-z^{n+b}|.
		\end{align*}
		On the other hand, using Lemma \ref{lem:2} we have 
		\begin{align*}
			d(Tx,Ty)&=|z^{n}-z^{n+1}|=2^n|1-z|=2^n\sqrt{7}=\frac{2}{\sqrt7}2^{n-1}\cdot7\\
			&\le \frac{2}{\sqrt7}|z^{n-1}||z^{a+1}-z^{b}|= \frac{2}{\sqrt7}|z^{n+a}-z^{n-1+b}|.
		\end{align*}
		The above two cases conclude the proof.
	\end{itemize}	
	\end{proof}
	\begin{remark}
		The proofs of Lemma \ref{lem:1}, \ref{lem:2} and Claim \ref{cl:2} can be carried out with fewer steps than presented (some cases can be merged). However, we think that these shortenings detriment the readability of the paper.
	\end{remark}
	
	\bibliographystyle{plain}
	\bibliography{references}

\end{document}